\documentclass{amsart}
\usepackage{graphicx}
\newtheorem{thm}{Theorem}[section]
\newtheorem{cor}[thm]{Corollary}

\theoremstyle{definition}

\theoremstyle{remark}

\numberwithin{equation}{section}

\newcommand{\R}{\mathbb R}
\newcommand{\Z}{\mathbb Z}
\newcommand{\ad}{\operatorname{ad}}

\newcommand{\range}{\operatorname{range}}
\newcommand{\Lie}{\mathfrak}

\begin{document}

\title[Graded Levi Decomposition]{The Levi Decomposition of a Graded Lie Algebra}%
\author{Paolo Ciatti}%
\address{Dipartimento di Ingegneria Civile, Edile e Ambientale, via Mazzolo 9, 35131 Padova, Italia}%
\email{paolo.ciatti@unipd.it}%

\author{Michael G. Cowling}
\address{School of Mathematics and Statistics, University of New South Wales, UNSW Sydney 2015, Australia}%
\email{m.cowling@unsw.edu.au}

\thanks{This work was developed while the second author was supported by a Visiting Scientist Grant from the Universit\`a di Padova}%

\begin{abstract}
We show that a graded Lie algebra admits a Levi decomposition that is compatible with the grading.
\end{abstract}
\maketitle
\section{introduction}
We assume throughout that Lie algebras are real and finite-dimensional.
In 1905, E. E. Levi \cite{Levi} showed that every Lie algebra may be decomposed as a direct sum:
\[
\Lie{g} = \Lie{l} \oplus \Lie{r},
\]
where $\Lie{r}$ is the radical, the maximal solvable ideal of $\Lie{g}$, and  $\Lie{l}$ is a semisimple subalgebra of $\Lie{g}$.
This is one of the foundations of Lie theory.

In more recent times, $\Z$-graded Lie algebras have come to assume an important role: for such an algebra $\Lie{g}$, we may write $\Lie{g} = \sum_{n\in \Z} \Lie{g}_n$, where $\Lie{g}_n = \{0\}$ for all but finitely many $n$ and $[\Lie{g}_m, \Lie{g}_n] \subseteq \Lie{g}_{m+n}$ for all $m$ and $n$.
Given a grading, there is an associated derivation, $\delta_1$ say, which is determined by linearity and the condition that $\delta_1 X = nX$ for all $X \in \Lie{g}_n$.
Conversely, given a diagonalisable derivation $\delta_1$, all of whose eigenvalues are integers, we obtain a grading by
defining $\Lie{g}_n$ to be the eigenspace for the eigenvalue $n$ when $n$ is an eigenvalue, and $\{0\}$ otherwise.

To he best of our knowledge, the interplay between the grading of an algebra and the Levi decomposition has not been made explicit, and this paper fills this gap.
Suppose that $\Lie{g}$ is a $\Z$-graded Lie algebra.
The radical $\Lie{r}$ is a characteristic ideal of $\Lie{g}$, so $\delta_1 \Lie{r} \subseteq \Lie{r}$, hence $\Lie{r} = \sum_n \Lie{r}_n$, where $\Lie{r}_n = \Lie{g}_n \cap \Lie{r}$.
But the same need not hold for a generic Levi subalgebra $\Lie{l}$.
We will show how to choose a $\delta_1$-invariant Levi subalgebra $\Lie{l}$; for this choice of $\Lie{l}$ it follows immediately that $\Lie{l} = \sum_n \Lie{l}_n$, where $\Lie{l}_n = \Lie{g}_n \cap \Lie{l}$.

More generally, Lie algebras may be graded over $\Z^d$, where $d > 1$; root space decompositions of semisimple Lie algebras are examples of this.
These gradings correspond to commuting families of diagonalisable derivations with integer eigenvalues.

We consider a slightly more general structure: an abelian Lie algebra $\Lie{a}$ and a homomorphism $\delta: H \mapsto \delta_H$ from $\Lie{a}$ into the Lie algebra of derivations of the Lie algebra $\Lie{g}$; we assume that each $\delta_H$ is semisimple, but not necessarily diagonalisable over $\R$.
Here is our main result.

\begin{thm}\label{thm-1}
Let $\{ \delta_H : H \in \Lie{a} \}$ be a commuting family of semisimple derivations of a Lie algebra $\Lie{g}$, and let $\Lie{r}$ be the radical of $\Lie{g}$.
Then there exists a Levi subalgebra $\Lie{l}$ such that $\delta_{\Lie{a}} \Lie{l} \subseteq \Lie{l}$.
\end{thm}

There are a number of arguments in the literature that this theorem illuminates and simplifies: see, for example, \cite{Cowling-Ottazzi} and \cite{Medori-Nacinovich}.
We add that Alexey Gordienko kindly pointed out to us that this result can be found in his papers \cite{Go1, Go2}, which treat Hopf algebraic questions; however our approach is more direct and uses Lie theory only.

\section{{Proof of Theorem~\ref{thm-1}}}
We follow the standard proof of the Levi decomposition, all the while keeping an eye on the derivations $\delta_H$.
First of all, we suppose that the derivations are inner, that is, $\Lie{a} \subset \Lie{g}$ and $\delta_H = \ad(H)$ for each $H \in \Lie{a}$.
Then $[\Lie{g}, \Lie{r}]$ and $[\Lie{r}, \Lie{r}]$ are $\ad(\Lie{a})$-invariant ideals in $\Lie{g}$ that are contained in $\Lie{r}$.
We consider several cases.

\textbf{Case 1}: There is an $\ad(\Lie{a})$-invariant ideal $\Lie{i}$ such that $ \{0\} \subset \Lie{i} \subset \Lie{r}$.

In this case, we argue by induction on dimension.
The derivations $\delta_H$ induce derivations on $\Lie{g} / \Lie{i}$, and $\Lie{r} / \Lie{i}$ is the radical of $\Lie{g} / \Lie{i}$, so we may write $\Lie{g} / \Lie{i} = \Lie{h} / \Lie{i} \oplus \Lie{r} / \Lie{i}$, where $\Lie{h}$ contains $\Lie{i}$ and is $\ad(\Lie{a})$-invariant, and $\Lie{h} / \Lie{i}$ is semisimple.
Then $\Lie{i}$ is the radical of $\Lie{h}$, and by induction, we may write $\Lie{h}$ as $\Lie{l} \oplus \Lie{i}$, where $\Lie{l}$ is semisimple and $\ad(\Lie{a})$-invariant.
Then $\Lie{g} = \Lie{l} \oplus \Lie{r}$ is an $\ad(\Lie{a})$-invariant Levi decomposition of $\Lie{g}$, and the result is established in Case 1.

Since $[\Lie{g}, \Lie{r}]$ is an $\ad(\Lie{a})$-invariant ideal in $\Lie{g}$, if we are not in Case 1, then either $[\Lie{g}, \Lie{r}] = \{0\}$ or  $[\Lie{g}, \Lie{r}] = \Lie{r}$.
Similarly, either $[\Lie{r}, \Lie{r}] = \{0\}$ or  $[\Lie{r}, \Lie{r}] = \Lie{r}$; this latter case cannot occur as $\Lie{r}$ is solvable.

\textbf{Case 2a}: $[\Lie{g}, \Lie{r}] = \{0\}$.

In this case, $\Lie{r}$ is the centre of $\Lie{g}$, and $\Lie{g} = \Lie{l} \oplus \Lie{r}$, where $\Lie{l} = [\Lie{g}, \Lie{g}]$; both summands are $\ad(\Lie{a})$-invariant and $\Lie{l}$ is a Levi subalgebra.

\textbf{Case 2b}: $[\Lie{g}, \Lie{r}] = \Lie{r}$.

In this case, $[\Lie{r}, \Lie{r}] = \{ 0 \}$, so $\ad(X)^2 = 0$ for all $X \in \Lie{r}$; moreover, the centre of $\Lie{g}$ is trivial, for otherwise we are in case 1.
We take a Levi decomposition $\Lie{l} \oplus \Lie{r}$ of $\Lie{g}$, and modify $\Lie{l}$ to achieve the desired decomposition.

Take $H \in \Lie{a} \setminus \{0\}$, and write $H = H_{\Lie{l}} + H_{\Lie{r}}$, where $H_{\Lie{l}} \in \Lie{l}$ and $H_{\Lie{r}} \in \Lie{r}$.
Since $\ad(H)$ maps $\Lie{r}$ into $\Lie{r}$ and  is semisimple, $\Lie{r} = \ker \ad(H) \oplus \range \ad(H)$,
so we may write $H_{\Lie{r}}$ as $H_0 + \ad(H) X$, where $H_0 \in \ker(\ad(H))$ and $X \in \Lie{r}$.

Now $\exp(\ad(X))$ is an automorphism of $\Lie{g}$, so we may define a new Levi factor, $\tilde{\Lie{l}}$ say, to be $\exp(\ad(X)) \Lie{l}$; we also define $\tilde H$ and $\tilde H_{\Lie{l}}$ to be $\exp(\ad(X)) H$ and $\exp(\ad(X)) H_\Lie{l}$.
Now $\ad(X) H_{\Lie{l}} = - [H_{\Lie{l}}, X] = - [H, X]$ as $[H_{\Lie{r}}, X] = 0$, whence
\[
\tilde H_{\Lie{l}} = H_{\Lie{l}} + \ad(X) H_{\Lie{l}} = H_{\Lie{l}} - [H, X] = H - H_{\Lie{r}} - [H,X] = H - H_0.
\]
Hence $H = \tilde H_{\Lie{l}} + H_0$.

By definition, $[H, H_0] = 0$, so $[H, \tilde H_{\Lie{l}}] = 0$.
Moreover, $\ad(H)$ is semisimple by definition.
Further, the action of $\ad(\tilde H_{\Lie{l}})$ on $\tilde{\Lie{l}}$ coincides with the quotient action of $\ad(\tilde H)$ on $\tilde{\Lie{l}} \oplus\Lie{r} / \Lie{r}$, which is semisimple by definition, whence the action of $\ad(\tilde H_{\Lie{l}})$ is also semisimple on $\Lie{g}$ (see, for example, \cite[Corollary C.18, p.~483]{Fulton-Harris}).
As $\ad(H)$ and $\ad(\tilde H_{\Lie{l}})$ commute and are both semisimple, $\ad(H - \tilde H_{\Lie{l}})$ is also semisimple.
However, $\ad(H - \tilde H_{\Lie{l}})= \ad(H_0)$ and $\ad(H_0)$ is nilpotent as $H_0 \in \Lie{r}$.
We deduce that $\ad(H_0)$ is trivial.
As the centre of $\Lie{g}$ is trivial, $H_0 = 0$.
In conclusion, $H = \tilde H_{\Lie{l}} \in \tilde{\Lie{l}}$, and $\tilde{\Lie{l}}$ is $\ad(H)$-invariant.

This argument shows that we can take a Levi subalgebra $\Lie{l}$ that is $\ad(H)$-invariant for a given $H$ in $\Lie{a}$, but $\Lie{a}$ may not be $1$-dimensional, so more is required.
After passing to the complexification if necessary, we may suppose that $\Lie{r} = \sum_{\alpha \in \Sigma} \Lie{r}_\alpha$, where $\Sigma$ is a finite subset of $\Lie{a}^*$, and $[H, X] = \alpha(H) X$ for all $X \in \Lie{r}_\alpha$ and all $H \in \Lie{a}$.
We take $H \in \Lie{a}$ such that $\alpha(H) \neq 0$ for all $\alpha \in \Sigma$, and assume that $H \in \Lie{l}$ so that $\Lie{l}$ is $\ad(H)$-invariant.

Now if $H' \in \Lie{a}$, and we write $H' = H'_{\Lie{l}} + H'_{\Lie{r}}$, where $H'_{\Lie{l}} \in \Lie{l}$ and $ H'_{\Lie{r}} \in \Lie{r}$, then $H$ and $H'$ commute by definition, and $H_{\Lie{r}}$ and $H'_{\Lie{r}}$ commute because $\Lie{r}$ is abelian.
Further, using the identification of $\Lie{l}$ with $\Lie{g} / \Lie{r}$, we see that $\ad(H)$ and $\ad(H')$ induce commuting derivations of $\Lie{l}$, which may be identified with $\ad(H_{\Lie{l}})$ and $\ad(H'_{\Lie{l}})$, whence $H_{\Lie{l}}$ and $H'_{\Lie{l}}$ commute.
Write $H_{\Lie{r}}$ as $\sum_\alpha H_\alpha$ and $H'_{\Lie{r}}$ as $\sum_\alpha H'_\alpha$, where $H_\alpha, H'_\alpha \in \Lie{r}_\alpha$.
Then
\[
\begin{aligned}
0
= [H, H']
&= [H_{\Lie{l}} + H_{\Lie{r}}, H'_{\Lie{l}} + H'_{\Lie{r}}]
= [H_{\Lie{l}}, H'_{\Lie{r}}] + [H_{\Lie{r}}, H'_{\Lie{l}}] \\
&= [H , H'_{\Lie{r}}] - [H', H_{\Lie{r}}]
= \sum_\alpha \alpha(H) H'_{\alpha} - \sum_\alpha \alpha(H') H_{\alpha}.
\end{aligned}
\]
Since the ``root space'' decomposition of $\Lie{r}$ is a direct sum, $\alpha(H) H'_{\alpha} = \alpha(H') H_\alpha$.
But $H_\alpha = 0$ for all $\alpha$, since $H \in \Lie{l}$, whence $H'_\alpha = 0$ unless $\alpha = 0$.
Finally, much as argued above, $\ad(H'_0)$ is both semisimple and nilpotent, hence null, whence $H'_0 = 0$, and $H' \in \Lie{l}$, as required.

We conclude our discussion of Cases 1 to 2b by affirming that Theorem 1 holds when all the derivations $\delta_H$ are inner.

It remains to discuss the general case, where some or all of the derivations are not inner.
In this case, we define ${\Lie{g}}_1$ to be the vector space $\Lie{g} \oplus \Lie{a}$, and consider $\Lie{g}$ and $\Lie{a}$ as subspaces of $\Lie{g}_1$ in the usual way.
We take the Lie product $[ \cdot, \cdot ]_1$ on $\Lie{g}_1$ that is determined by linearity, antisymmetry, and the requirements that
\[
[H,H'] = 0, \quad [X,X']_1 = [X,X'] \quad\text{and}\quad [H,X]_1 = \delta_H X
\]
for all $H,H' \in \Lie{a}$ and all $X, X' \in \Lie{g}$.

Now there exist a semisimple subalgebra $\Lie{l}_1$ and a solvable ideal $\Lie{r}_1$, both $\ad(\Lie{a})$-invariant, such that $\Lie{g}_1 = \Lie{l}_1 \oplus \Lie{r}_1$.
By construction, $\Lie{l}_1 = [\Lie{l}_1, \Lie{l}_1] \subseteq [\Lie{g}_1, \Lie{g}_1] \subseteq \Lie{g}$, so
$\Lie{l}_1$ is a $\delta_{\Lie{a}}$-invariant semisimple subalgebra of $\Lie{g}$.
Further, it is easy to see that $\Lie{g} \cap \Lie{r}_1 = \Lie{r}$.
Finally, since $\Lie{g}_1 = \Lie{l}_1 \oplus \Lie{r}_1$ and $\Lie{l}_1 \subseteq \Lie{g}$, it follows that
\[
\Lie{g} = \Lie{l}_1 \oplus \bigl(\Lie{g} \cap \Lie{r}_1\bigr) = \Lie{l}_1 \oplus \Lie{r},
\]
where both summands are $\delta_{\Lie{a}}$-invariant, as required.

\section{Some corollaries}

The first corollary is immediate.

\begin{cor}
Suppose that $\Lie{g}$ is a $\Z^d$-graded Lie algebra, so that $\Lie{g} = \sum_m \Lie{g}_m$, where $[\Lie{g}_m, \Lie{g}_n] \subseteq \Lie{g}_{m+n}$ for all $m, n  \in \Z^d$.
Then $\Lie{r} = \sum_m \Lie{r}_m$, where $\Lie{r}_m = \Lie{r} \cap \Lie{g}_m$, and there is a Levi subalgebra $\Lie{l}$ such that $\Lie{l} = \sum_m \Lie{l}_m$, where $\Lie{l}_m = \Lie{l} \cap \Lie{g}_m$.
\end{cor}

Our decomposition also provides us with some information about algebras of derivations.

\begin{cor}
Suppose that $\delta_{\Lie{a}}$ is an abelian algebra of semisimple derivations of $\Lie{g}$, and that $\Lie{l} \oplus \Lie{r}$ is a Levi decomposition of $\Lie{g}$ into $\delta_{\Lie{a}}$-invariant summands.
Then there are commuting algebras $\delta_{\Lie{a},\Lie{l}}$ and $\delta_{\Lie{a},\Lie{r}}$ of commuting semisimple derivations of $\Lie{g}$ that preserve $\Lie{l}$ and $\Lie{r}$ such that $\delta_{\Lie{a}} \subseteq \delta_{\Lie{a},\Lie{l}} \oplus \delta_{\Lie{a},\Lie{r}}$; further, $\delta_{\Lie{a},\Lie{l}}$ may be identified with $\ad(\Lie{b})$, where $\Lie{b}$ is a commutative subalgebra of $\Lie{l}$, and every element of $\delta_{\Lie{a}, \Lie{r}}$ annihilates $\Lie{l}$.
\end{cor}

To see this, recall that each derivation $\delta_H$ of $\Lie{g}$ induces a derivation of $\Lie{l} \oplus \Lie{r} / \Lie{r}$, which we may identify with an inner derivation $\ad(H_\Lie{l})$ of $\Lie{l}$; take $\delta_{\Lie{a}, \Lie{l}}$ to be the algebra of all derivations of $\Lie{g}$ of the form $\ad(H_{\Lie{l}})$; these commute by construction.
Define $\delta_{H, \Lie{r}}$ to be the derivation $\delta_H - \ad(H_\Lie{l})$; then $\delta_{H, \Lie{r}} (X) = 0$ for all $X \in \Lie{l}$.
It follows that $[\delta_{H, \Lie{r}}, \ad(X) ] = 0$ for all $X \in \Lie{l}$, and so
\[
 [ \delta_{H, \Lie{r}},  \delta_{H', \Lie{r}}] = [\ad(H_\Lie{l}) + \delta_{H, \Lie{r}}, \ad(H'_\Lie{l}) + \delta_{H', \Lie{r}}] = [\delta_H, \delta_{H'}]  = 0.
\]
In conclusion, the various $\delta_{H, \Lie{r}}$ commute amongst themselves and with the $\ad(H_{\Lie{l}})$, which also commute amongst themselves.
The corollary follows.

This corollary has implications for gradings of general Lie algebras: they arise from gradings on a Levi subalgebra and from gradings of the radical which are invariant under the action of the Levi subalgebra.

Finally, if $D$ is a derivation of a Lie algebra $\Lie{g}$, and annihilates the Levi subalgebra $\Lie{l}$, then $D$ is supported on the nilradical $\Lie{n}$ of $\Lie{G}$.
Indeed, $\ad(\Lie{l}) \oplus \R D$ is a reductive algebra of derivations of $\Lie{g}$, which stabilises $\Lie{r}$ and $\Lie{n}$, so there is an $\ad(\Lie{l}) \oplus \R D$-invariant subspace $\Lie{a}$ such that $\Lie{r} = \Lie{a} \oplus \Lie{n}$.
Further, $D \Lie{r} \subseteq \Lie{n}$ (see \cite[Corollary C.24, p.~485]{Fulton-Harris}), and so $D\ad{\Lie{a}} \subseteq \Lie{a} \cap \Lie{n}$, that is, $D|_{\Lie{a}} = 0$.


\section{Afterword}
Ironically, although the Levi decomposition appears in many textbooks and research papers, at the time of writing of this article, Levi's paper has been cited 5 times, according to \emph{Zentralblatt f\"ur Mathematik}; so much for the existence of a correlation between citation numbers and significance of a contribution.

\bibliographystyle{amsplain}

\end{document}